\input amstex
\input amsppt.sty
\magnification=\magstep1
\hsize=30truecc
\baselineskip=16truept
\vsize=22.2truecm
\NoBlackBoxes
\nologo
\pageno=1
\topmatter
\TagsOnRight

\def\N{\Bbb N}
\def\Z{\Bbb Z}

\def\l{\left}
\def\r{\right}
\def\b{\bigg}

\def\({\b(}
\def\[{\b[}
\def\){\b)}
\def\]{\b]}

\def\t{\text}
\def\f{\frac}
\def\mo{\roman{mod}}

\def\em{\emptyset}
\def\se {\subseteq}

\def\eq{\equiv}

\def\ls{\leqslant}
\def\gs{\geqslant}
\def\al{\alpha}

\def\da{\delta}

\def\Proof{\noindent{\it Proof}}
\def\Remark{\noindent{\it Remark}}

\hbox {Journal of Combinatorics and Number Theory {\bf 1}(2009),
no.\,1, 65--76.}
\medskip
\title On sums of primes and triangular numbers \endtitle
\author  Zhi-Wei Sun \endauthor

\affil Department of Mathematics, Nanjing University
     \\Nanjing 210093, People's Republic of China
    \\ {\tt zwsun\@nju.edu.cn}
    \\ http://math.nju.edu.cn/$\sim$zwsun
 \endaffil
\abstract We study whether sufficiently large integers can be written in the form $cp+T_x$, where
$p$ is either zero or a prime congruent to $r$ mod $d$, and
$T_x=x(x+1)/2$ is a triangular number. We also investigate whether there are infinitely many
positive integers not of the form $(2^ap-r)/m+T_x$ with $p$ a prime and $x$ an integer.
Besides two theorems, the paper also contains several conjectures
together with related analysis and numerical data. One of our conjectures states that
each natural number $n\not=216$ can be written in the form $p+T_x$ with $x\in\Z$ and $p$ a prime or zero;
another conjecture asserts that any odd integer $n>3$ can be written in the form $p+x(x+1)$
with $p$ a prime and $x$ a positive integer.
\endabstract
\thanks Received on April 9, 2008, and accepted on June 19, 2008.
\newline\indent Supported by the National Natural Science Foundation (grant 10871087) of People's Republic of China.
\endthanks
\endtopmatter
\document

\heading{1. Introduction}\endheading

 Since $1+2+\cdots+n=n(n+1)/2$, those integers $T_x=x(x+1)/2$ with
 $x\in\Z$ are called {\it triangular numbers}. Note that $T_{-x}=T_{x-1}$.
 Here is a list of the triangular numbers
not exceeding 200:
$$0,\,1,\, 3,\, 6,\, 10,\, 15,\, 21,\, 28,\ 36,\, 45,\, 55,\,
66,\, 78,\, 91,\, 105,\, 120,\, 136,\, 153,\, 171,\, 190.$$
For $n,x\in\N=\{0,1,2,\ldots\}$, clearly
$$T_n\ls x\iff 2n+1\ls\sqrt{8x+1}.$$
Thus, for any $x\gs 0$, there are exactly $\lfloor (\sqrt{8x+1}-1)/2\rfloor+1$
triangular numbers not exceeding $x$.

 Here is an important observation of Fermat.

\proclaim{Fermat's Assertion} Each $n\in\N$ can be written as a sum of three triangular numbers.
\endproclaim

 An equivalent version of this assertion states that for each $n\in\N$ the number
 $8n+3$ is a sum of three squares (of odd integers). This is a consequence of
 the following profound theorem (see, e.g.,
 [Gr, pp.\,38--49] or [N, pp.\,17-23]) due to Gauss and Legendre:
 A natural number can be written as a sum of three squares
 of integers if and only if it is not of the form $4^k(8l+7)$ with $k,l\in\N$.

 Prime numbers play a key role in number theory.
 By the prime number theorem, for $x\gs 2$ the number $\pi(x)$ of primes not exceeding $x$
 is approximately $x/\log x$ (in fact, $\lim_{x\to+\infty}\pi(x)/(x/\log x)=1$).
 Here is a famous result due to I. M. Vinogradov [V].

\proclaim{Vinogradov's Theorem} Every sufficiently large odd integer can be written as a sum of three primes.
\endproclaim

  The following result of Linnik [L1, L2] is also remarkable: Any sufficiently large integer can be written
 as a sum of a prime and two squares of integers.

 Now we state a well-known conjecture which remains unsolved.

 \proclaim{Goldbach's Conjecture} Any even number greater than two can be expressed as a sum of two primes.
 \endproclaim

 In this paper we investigate mixed sums of primes and triangular numbers. It seems that no one has studied this topic before.
 Surprisingly, there are many mysterious things in this new field.

 Here is our first result.
\proclaim{Theorem 1.1} Let $c,d\in\Z^+=\{1,2,3,\ldots\}$ and $r\in\Z$. Assume that there are only finitely many natural numbers
not in the form $cp+T_x$, where $p$ is zero or a prime in the residue class $r(\mo\ d)$, and $x$ is an integer.
Then both $c$ and $d$ are powers of two, and $r$ is relatively prime to $d$.
\endproclaim

For any $d\in\Z^+$ and $r\in\Z$ with $(r,d)=1$ (where $(r,d)$ denotes the greatest common divisor of $r$ and $d$),
the residue class $r(\mo\ d)$ contains infinitely many primes by Dirichlet's theorem (cf. [IR, p.\,251]).

Let $a\in\N$ and $r\in\Z$. It is known that $8r+1$ is a quadratic
residue modulo $2^{a+3}$ (see, e.g., [IR, Proposition 5.1.1]). So,
there is an integer $x$ such that $(2x+1)^2\eq 8r+1\ (\mo\
2^{a+3})$, i.e., $T_x\eq r\ (\mo\ 2^a)$. Thus $\{T_x:\ x\in\Z\}$
contains a complete system of residues modulo any power of two.
This important property of triangular numbers leads us to make the
following conjecture in view of Theorem 1.1.

\proclaim{Conjecture 1.1} For any $a,b\in\N$ and odd integer $r$, all sufficiently large integers
can be written in the form $2^ap+T_x$ with $x\in\Z$, where $p$ is either zero or a prime congruent to $r$ mod $2^b$.
In particular, each natural
number $n\not=216$ can be written in the form $p+T_x$ with $x\in\Z$, where $p$
is zero or a prime; furthermore, any positive integer $n\not\in\{2,5,7,61,211,216\}$ can be written
in the form $p+T_x$ with $x\in\Z^+$, where $p$ is an odd prime or zero.
\endproclaim

\Remark\ 1.1. (i) Conjecture 1.1 seems quite unexpected,
nevertheless we have verified its latter part for $n\ls 17,000,000$.
It is interesting to compare Conjecture 1.1 with the Goldbach
conjecture. Note that there are much more primes than triangular
numbers below large $x$. As for the number 216, it is well known
that $216=6^3=3^3+4^3+5^3$. (ii) In March and April, 2008, the
author posted several messages concerning Conjecture 1.1 and related
things to the Number Theory Mailing List; the first of which was
made public on March 23, 2008 (cf. {\tt
http://listserv.nodak.edu/cgi-bin/wa.exe?A2=ind0803\&L=nmbrthry\&T=0\&P=3010}).
\medskip

A well-known assertion of Fermat (proved by Euler) states that
each prime $p\eq1\ (\mo\ 4)$ can be written in the form $x^2+y^2$ with $x$ even and $y$ odd
(see, e.g., [G, pp.\,163--165] or [IR, p.\,64]). Thus Conjecture 1.1 implies that
for any $a=0,1,2,\ldots$ all sufficiently large integers have the form $2^a(x^2+y^2)+T_z$ with $x,y,z\in\Z$.
It is known that, if a positive integer is not a triangular number,
 then it must be a sum of an even square, an odd square and a triangular number (cf. [S, Theorem 1(iii)]).
If $p=x^2+y^2$ with $x$ even and $y$ odd, then $2p=(x+y)^2+(x-y)^2$ with $x\pm y$ odd.
Thus our following conjecture is reasonable in view of Conjecture 1.1.

\proclaim{Conjecture 1.2}  {\rm (i)} A natural number can be written as a sum of two even squares and a triangular number
unless it is among the following list of $19$ exceptions:
$$\align &2,\ 12,\ 13,\, 24,\ 27,\ 34,\ 54,\ 84,\ 112,\ 133,
\\&162,\ 234,\ 237,\ 279,\ 342,\  399,\ 652,\  834,\ 864.\endalign$$
Furthermore, any integer $n>2577$ can be written in the form $(4x)^2+(2y)^2+T_z$ with $x,y,z\in\Z$.

{\rm (ii)} Each natural number $n\not\in E$ is either a triangular number, or a sum of a triangular number and two odd squares,
where the exceptional set $E$ consists of the following $25$ numbers:
$$\align &4,\ 7,\ 9,\ 14,\ 22,\ 42,\ 43,\ 48,\ 52,\ 67,\ 69,\ 72,\ 87,\ 114,
\\&144,\ 157,\ 159,\ 169,\ 357,\ 402,\ 489,\ 507,\ 939,\ 952,\ 1029.
\endalign$$
\endproclaim

\Remark\ 1.2. We have verified Conjecture 1.2 for $n\ls 2,000,000$.
In [S] the author proved that any natural number $n$ is a sum of an even square and two triangular numbers,
and we also can express $n$ as a sum of an odd square and two triangular numbers if it is not twice a triangular number.
For other problems and results on mixed sums of squares and triangular numbers, the reader may consult [S] and [GPS]
and the references therein.

\medskip

Let $m>1$ be an integer, and let $a\in\N$ with $(2^a,m)=1$.
We define
 $$S^{(a)}_m=\{n>m:\ (m,n)=1\ \&\ n\not=2^ap+mT_x\ \t{for any prime}\ p\ \t{and}\ \t{integer}\ x\}$$
 and simply write $S_m$ for $S_m^{(0)}$.
Clearly
$$S^{(a)}_m=\bigcup\Sb 1\ls r\ls m\\(r,m)=1\endSb \{r+mn:\ n\in S^{(a)}_m(r)\},$$
where
$$S^{(a)}_m(r)=\l\{n\in\Z^+: \ n\not=\f{2^ap-r}m+T_x\ \t{for any prime}\ p\ \t{and integer}\ x\r\}.$$
(We also abbreviate $S^{(0)}_m(r)$ to $S_m(r)$.)
What can we say about these exceptional sets? Are they finite?

Here is our second theorem.

\proclaim{Theorem 1.2} {\rm (i)} Let $m>1$ be an odd integer, and let $a\in\N$.
If $r$ is a positive integer such that
$2r$ is a quadratic residue modulo $m$, then there are infinitely many positive integers
not of the form $(2^ap-r)/m^2+T_x$, where $p$ is a prime and $x$ is an integer.
Therefore the set $S^{(a)}_{m^2}$ is infinite.

{\rm (ii)} Let $m=2^{\al}m_0$ be a positive even integer with $\al,m_0\in\Z^+$ and $2\nmid m_0$.
 If $r\in\Z^+$ is a quadratic residue modulo $m_0$ with
$r\eq 2^{\al}+1\ (\mo\ 2^{\min\{\al+1,3\}})$, then there are infinitely many
positive integers not of the form $(p-r)/(2m^2)+T_x$, where $p$ is a prime and $x$ is an integer. Thus
$S_{2m^2}$ is an infinite set.
\endproclaim

\Remark\ 1.3. Let $m$ be a positive odd integer. By a well known result (see, e.g., [IR, pp.\,50-51]),
an integer $r$ is a quadratic residue mod $m$ if and only if for any prime divisor $p$ of $m$ the Legendre symbol
$(\f rp)$ equals one.

\medskip

In view of Theorem 1.2 and some computational results, we raise the following conjecture.

\proclaim{Conjecture 1.3} Let $m>1$ be an integer.

 {\rm (i)} Assume that $m$ is odd. If $m$ is not a square, then $S_m,\, S_m^{(1)},\ S_m^{(2)},\ldots$ are all finite.
If $m=m_0^2$ with $m_0\in\Z^+$, and $r$ is a positive integer with $(r,m)=1$
such that $2r$ is a quadratic non-residue mod $m_0$,
then $S_m^{(a)}(r)$ is finite for every $a=0,1,2,\ldots$.

 {\rm (ii)} Suppose that $m$ is even. If $m$ is not twice an even square, then the set $S_m$ is finite.
If $m=2(2^{\alpha}m_0)^2$ with $\al,m_0\in\Z^+$ and $2\nmid m_0$, and $r$ is a positive integer with $(r,m)=1$
 such that $r$ is a quadratic non-residue modulo $m_0$ or $r\not\eq 2^{\al}+1\ (\mo\ 2^{\min\{\alpha+1,3\}})$,
 then the set $S_m(r)$ is finite.
\endproclaim

\noindent {\it Example} 1.1. (i) Among $1,\ldots,15$ only $1$ and $4$ are quadratic residues modulo 15.
For any $a\in\N$, both $S_{15^2}^{(a)}(2)$ and $S_{15^2}^{(a)}(8)$ are infinite by Theorem 1.2(i), while
$$S_{225}^{(a)}(1),\ S_{225}^{(a)}(4),\ S_{225}^{(a)}(7),\ S_{225}^{(a)}(11),\ S_{225}^{(a)}(13),\ S_{225}^{(a)}(14)$$
should be finite as predicted by Conjecture 1.3(i).

(ii) Let $r$ be a positive odd integer. By Theorem 1.2(ii), $S_{2\times 8^2}(r)$ is infinite if $r\eq1\ (\mo\ 8)$.
On the other hand, by Conjecture 1.3, $S_{2\times 8^2}(r)$ should be finite when $r\not\eq1\ (\mo\ 8)$.
When $(r,18)=1$, the set $S_{2\times18^2}(r)$ is infinite if $r\eq 7\ (\mo\ 12)$ (i.e., $r$ is quadratic residue mod $3^2$
with $r\eq3\ (\mo\ 4)$) (by Theorem 1.2(ii)),  and it is finite otherwise (by Conjecture 1.3(ii)).
Similarly, when $(r,20)=1$, the set $S_{2\times20^2}(r)$ is infinite if $r\eq 21,29\ (\mo\ 40)$ (i.e., $r$ is quadratic residue mod $5$
with $r\eq5\ (\mo\ 8)$) (by Theorem 1.2(ii)),  and it is finite otherwise (by Conjecture 1.3(ii)).

\bigskip

In view of Conjecture 1.3(ii), the sets $S_2,S_6,S_{12}$ and $S_{288}(19)$ should be finite.
Our computations up to $10^6$ suggest further that $S_2=S_6=S_{12}=S_{288}(19)=\em$.
Recall that
$$S_2=\{2n+1:\ n\in\Z^+\ \t{and}\ 2n+1\not=p+2T_x\ \t{for any prime}\ p\ \t{and integer}\ x\}.$$

 Now we pose one more conjecture.
\proclaim{Conjecture 1.4} Any odd integer $n>3$ can be written in the form $p+x(x+1)$
with $p$ a prime and and $x$ a positive integer. Furthermore, for any $b\in\N$ and $r\in\{1,3,5,\ldots\}$
all sufficiently large odd integers can be written in the form $p+x(x+1)$ with $x\in\Z$,
where $p$ is a prime congruent to $r$ mod $2^b$.
\endproclaim

\Remark\ 1.4. It is interesting to compare the above conjecture with a conjecture of E. Lemoine (cf. [KY])
posed in 1894 which states that any odd integer greater than 5 can be written in  the form $p+2q$
where $p$ and $q$ are primes.

 \medskip

 In the next section we are going to prove Theorems 1.1 and 1.2. Section 3 is devoted to
 numerical illustrations of Conjectures 1.1 and 1.4.
 In Section 4 we present some additional remarks on Conjectures 1.2 and 1.3.

 \heading{2. Proofs of Theorems 1.1 and 1.2}\endheading

\proclaim{Lemma 2.1} Let $p$ be an odd prime. Then
$$|\{T_x\ \mo\ p:\ x\in\Z\}|=\f{p+1}2.\tag2.1$$
\endproclaim
\Proof.  For any $r\in\Z$, clearly
$$\align& T_x\eq r\ (\mo\ p)\ \t{for some}\ x\in\Z
\\\iff&(2x+1)^2=8T_x+1\eq 8r+1\ (\mo\ p)\ \t{for some}\ x\in\Z
\\\iff& p\mid 8r+1\ \t{or}\ 8r+1\ \t{is a quadratic residue}\ \mo\ p.
\endalign$$ Therefore $|\{T_x\ \mo\ p:\ x\in\Z\}|=(p-1)/2+1=(p+1)/2$. \qed

 \medskip
 \noindent{\it Proof of Theorem 1.1}. If $(r,d)$ has a prime divisor $q$, then there is no prime $p\not=q$ in the
 residue class $r(\mo\ d)$, and hence there are infinitely many natural numbers not in the form
 $cp+T_x$ with $x\in\Z$ and $p\in\{0\}\cup\{\t{primes in}\ r(\mo\ d)\}\se\{0,q\}$,
 which contradicts the assumption in Theorem 1.1. Therefore we have $(r,d)=1$.

 Suppose that $cd$ has an odd prime divisor $q$.
 As $(q+1)/2<q$, by Lemma 2.1 there is an integer $y$ with $y\not\eq cr+T_x\ (\mo\ q)$
for any $x\in\Z$. For any $n\in\N$, if we can write $y+nq$ in the form $cp+T_x$ with $x\in\Z$,
where $p$ is zero or a prime congruent to $r$
mod $d$, then $p$ must be zero, for, otherwise $y-T_x\eq y+nq-T_x=cp\eq cr\ (\mo\ q)$ which is impossible by the choice of $y$.
As there are infinitely many positive integers in the residue class $y(\mo\ q)$ which are not triangular numbers,
we get a contradiction and this concludes the proof. \qed

\medskip
\noindent{\it Proof of Theorem 1.2}.
(i) Suppose that $r$ is a positive integer for which $2r$ is a quadratic residue modulo the odd integer $m>1$.
Then there is an odd number $x\in\Z^+$ such that
$x^2\eq 2r\ (\mo\ m)$. (Note that if $x$ is even then $x+m$ is odd.) Thus $2r=x^2+mq$ for some odd integer $q$.
As $(x,m)=1=(m,2^{a+1})$, by the Chinese Remainder Theorem, for some integer $b\gs |q|$ we have
both $bx\eq q\ (\mo\ m)$ and $2x+bm\eq-1\ (\mo\ 2^{a+1})$.
Note that $b$ is odd since $bm\eq-1\ (\mo\ 2)$.
For $k\in\N$ we set
$$b_k=b+2^{a+1}km\quad\t{and}\quad n_k=\f{b_kx-q}{2m}+\f{b_k^2-1}8\in\N.\tag2.2$$
Then
$$(8n_k+1)m^2+8r=\l(b_k^2+4\f{b_kx-q}m\r)m^2+4(x^2+mq)=(2x+b_km)^2.$$

For every $k=1,2,3,\ldots$ we have
$$n_k\gs\f{b_k^2-1}8\gs\f{(2m+1)^2-1}8=\f{m(m+1)}2>m.$$
If $n_1,n_2,\ldots$ all belong to $S^{(a)}_{m^2}(r)$, then
$S^{(a)}_{m^2}(r)$ is obviously infinite.

Below we suppose that $S^{(a)}_{m^2}(r)$ does not contain all those
$n_1,n_2,\ldots$. Let $k$ be any positive integer with $n_k\not\in
S^{(a)}_{m^2}(r)$, i.e., $n_k=(2^ap-r)/m^2+T_z$ for some prime $p$
and $z\in\N$. Then $8n_k+1=8(2^ap-r)/m^2+y^2$, where $y=2z+1$ is a
positive odd integer. Therefore
$$2^{a+3}p=(8n_k+1)m^2+8r-(my)^2=(2x+b_km)^2-(my)^2.$$
Note that both $2x+b_km$ and $my$ are odd. As
$2x+b_km>2b_k\gs2^{a+2}$, for some $i\in\{0,\ldots,a+1\}$ we have
$$2x+b_km+my=2^{i+1}p\quad \t{and}\quad 2x+b_km-my=2^{a+2-i}.$$
(Note that $2x+b_km+my=p=2$ is impossible.)
It follows that
$$2x+b_km=\f{2^{i+1}p+2^{a+2-i}}2=2^ip+2^{a+1-i}.$$
Since $2x+b_km$ is odd, we must have $i\in\{0,a+1\}$. So
$$2x+b_km\in\{p+2^{a+1},\ 2^{a+1}p+1\}.\tag2.3$$

{\it Case}\ 1. $a>0$.

In this case,
$$2x+b_km\eq 2x+bm\eq-1\not\eq 1\ (\mo\ 2^{a+1}).$$
So $2x+b_km=p+2^{a+1}$. For each $l=1,2,3,\ldots$, obviously
$$2x+b_{k+lp}m-2^{a+1}=2x+b_km-2^{a+1}+2^{a+1}lpm^2=p(1+2^{a+1}lm^2)$$
and hence it is not a prime number. Therefore all the infinitely many numbers
$$n_{k+p}<n_{k+2p}<n_{k+3p}<\cdots$$
belong to the exceptional set $S_{m^2}^{(a)}(r)$.

\medskip

{\it Case}\ 2. $a=0$.

In this case, $2x+b_km\in\{p+2,\,2p+1\}$.

Assume that $2x+b_km=p+2$. Then, for each $l\in\Z^+$, we have
$$\align 2x+b_{k+lp(p+1)}m=&2x+b_km+lp(p+1)2m^2=p+2+2lm^2p(p+1)
\\=&2+p(1+2lm^2(p+1))=1+(p+1)(1+2lm^2p)
\endalign$$
and hence $2x+b_{k+lp(p+1)}m$ is not of the form $p'+2$ or $2p'+1$ with $p'$ a prime.
Thus, all the infinitely many positive integers
$n_{k+lp(p+1)}\ (l=1,2,3,\ldots)$ must belong to $S_{m^2}(r)$.

Now suppose that $2x+b_jm-2$ is not a prime for any $j\in\Z^+$ with $n_j\not\in S_{m^2}(r)$.
Then $2x+b_km=2p+1$. For each $l=1,2,3,\ldots$, the number
$$2x+b_{k+lp}m=2x+b_km+2lpm^2=2p(1+lm^2)+1$$
is not of the form $2p'+1$ with $p'$ a prime. So
all the infinitely many positive integers
$n_{k+lp}\ (l=1,2,3,\ldots)$ lie in $S_{m^2}(r)$.

\medskip

(ii) Now we proceed to the second part of Theorem 1.2.
Suppose that $r\in\Z^+$ is a quadratic residue mod $m_0$ with $r\eq 2^{\al}+1\ (\mo\ 2^{\min\{\al+1,3\}})$.
Note that the congruence $x^2\eq r\ (\mo\ 2^\al)$ is solvable since $r\eq1\ (\mo\ 4)$ if $\al=2$, and $r\eq1\ (\mo\ 8)$
if $\al\gs3$. By the Chinese Remainder Theorem, there is an integer $x$
such that $x^2\eq r\ (\mo\ 2^\al m_0)$ with $0<x\ls m/2=2^{\al-1}m_0$.
Write $r=x^2+mq$ with $q\in\Z$. If $8\nmid m$ (i.e., $\al\ls 2$), then
$$mq=r-x^2\eq(2^\al+1)-1\ (\mo\ 2^{\al+1})$$
and hence $q$ is odd.

Define
$$\delta=\cases0&\t{if}\ q\eq1\ (\mo\ 2),\\ 1&\t{if}\ 8\mid m\ \t{and}\ q\eq m/4\ (\mo\ 4),
\\5&\t{if}\ 8\mid m\ \t{and}\ q\eq m/4+2\ (\mo\ 4).\endcases$$
As $(x,2m)=1$, there is an integer $b>|q|$ such that
$$bx\eq q+\f m4\da(1-x)\ (\mo\ 2m)$$
and hence
$$\l(b+\f m4\da\r)^2\eq\l(bx+\f m4\da x\r)^2
\eq \l(q+\f m4\da\r)^2\eq1-\da\ (\mo\ 8).$$
For $k\in\N$ we set
$$b_k=b+\f m4\da+2km \quad\ \t{and}\quad \ n_k=\f{b_k^2+\delta-1}8+\f{b_kx-q-\da m/4}{2m}.\tag2.4$$
Clearly
$$b_k^2\eq\l(b+\f m4\da\r)^2\eq1-\da\ (\mo\ 8)$$
and $$b_kx\eq\l(b+\f m4\da\r)x\eq q+\f m4\da\ (\mo\ 2m).$$
So we have $n_k\in\Z$. Observe that
$$(8n_k+1)\l(\f m2\r)^2+r=(b_k^2+\delta)\f{m^2}4+m\l(b_kx-q-\da\f m4\r)+mq+x^2=\l(\f m2b_k+x\r)^2.$$

For $k\in\Z^+$, obviously $b_k\gs b+2m>|q|+2m$ and hence
$$n_k\gs\f{(2m+1)^2-1}8>m.$$
If $n_1,n_2,\ldots$ all belong to $S_{2m^2}(r)$, then $S_{2m^2}(r)$
is infinite.

Below we assume that there is a positive integer $k$ such that $n_k\not\in S_{2m^2}(r)$,
i.e., $n_k=(p-r)/(2m^2)+T_z$ for some prime $p$ and $z\in\N$. Then
$8n_k+1=(p-r)/(m/2)^2+y^2$ with $y=2z+1\in\Z^+$. Thus
$$p=(8n_k+1)\l(\f m2\r)^2+r-\l(\f m2y\r)^2=\l(\f m2b_k+x\r)^2-\l(\f m2y\r)^2.$$
Since $p$ is a prime, this implies that
$$\f m2b_k+x-\f m2y=1\ \ \t{and}\ \ \f m2b_k+x+\f m2y=p.$$
As $0<x\ls m/2$ and $x\eq 1\ (\mo\ m/2)$, we must have $x=1$ and $b_k=y$.
Therefore $b_km+1=p$.

 For each $l=1,2,3,\ldots$, clearly
$$b_{k+lp}m+1=b_km+1+2m^2lp=p(1+2lm^2)$$
is not a prime. Thus, by the above, all the infinitely many positive integers $n_{k+pl}\ (l=1,2,3,\ldots)$
must belong to the set $S_{2m^2}(r)$.

In view of the above, we have completed the proof of Theorem 1.2. \qed

\heading{3. Numerical illustrations of Conjectures 1.1 and 1.4}\endheading

Concerning the particular case $a=0$ and $b\in\{2,3\}$ of Conjecture 1.1, we have a more concrete conjecture.

\proclaim{Conjecture 3.1} {\rm (i)} Each natural number $n>88956$ can be written in the form
$p+T_x$ with $x\in\Z^+$, where $p$ is either zero or a prime congruent to $1$ mod $4$.
Each natural number $n>90441$ can be written in the form
$p+T_x$ with $x\in\Z^+$, where $p$ is either zero or a prime congruent to $3$ mod $4$.

{\rm (ii)} For $r\in\{1,3,5,7\}$, we can write any integer $n>N_r$
in the form $p+T_x$ with $x\in\Z$, where $p$ is either zero or a prime congruent to $r$ mod $8$, and
$$N_1=1004160,\ \ N_3=1142625,\ \ N_5=779646,\ \ N_7=893250.$$
\endproclaim

\Remark\ 3.1. We have verified Conjecture 3.1 for $n\ls 5,000,000$.
Since any prime $p\eq1\ (\mo\ 8)$ can be written in the form
$x^2+2(2y)^2$ with $x,y\in\Z$ (cf. [G, pp.\,165--166]), and all
natural numbers not exceeding $1,004,160$ can be written in the form
$x^2+8y^2+T_z$ with $x,y,z\in\Z$, Conjecture 3.1(ii) with $r=1$
implies the following deep result of Jones and Pall [JP] obtained by
the theory of ternary quadratic forms: For each natural number $n$
there are $x,y,z\in\Z$ such that $n=x^2+8y^2+T_z$, i.e.,
$8n+1=2(2x)^2+(8y)^2+(2z+1)^2$.

\medskip

Here is a list of all natural numbers not exceeding 88,956 that cannot be written in the form
$p+T_x$ with $x\in\Z$, where $p$ is either 0 or a prime congruent to 1 mod 4.
\medskip

2, 4, 7, 9, 12, 22, 24, 25, 31, 46, 48, 70, 75, 80, 85, 87, 93, 121, 126, 135, 148, 162,
169, 186, 205, 211, 213, 216, 220, 222, 246, 255, 315, 331, 357, 375, 396, 420, 432, 441,
468, 573, 588, 615, 690, 717, 735, 738, 750, 796, 879, 924, 1029, 1038, 1080, 1155, 1158,
1161, 1323, 1351, 1440, 1533, 1566, 1620, 1836, 1851, 1863, 1965, 2073, 2118, 2376, 2430,
2691, 2761, 3156, 3171, 3501, 3726, 3765, 3900, 4047, 4311, 4525, 4605, 4840, 5085, 5481,
5943, 6006, 6196, 6210, 6471, 6810, 6831, 6840, 7455, 7500, 7836, 8016, 8316, 8655, 8715,
8991, 9801, 10098, 10563, 11181, 11616, 12165, 12265, 13071, 14448, 14913, 15333, 15795,
17085, 18123, 20376, 27846, 28161, 30045, 54141, 88956.

\medskip

Below is a list of all natural numbers not exceeding 90,441 that cannot be written in the form
$p+T_x$ with $x\in\Z$, where $p$ is either 0 or a prime congruent to 3 mod 4.
\medskip

2, 5, 16, 27, 30, 42, 54, 61, 63, 90, 96, 129, 144, 165, 204, 216, 225, 285, 288, 309, 333, 340,
345, 390, 405, 423, 426, 448, 462, 525, 540, 556, 624, 651, 705, 801, 813, 876, 945, 960, 1056,
1230, 1371, 1380, 1470, 1491, 1827, 2085, 2157, 2181, 2220, 2355, 2472, 2562, 2577, 2655, 2787,
2811, 2826, 2886, 3453, 3693, 3711, 3735, 3771, 3840, 3981, 4161, 4206, 4455, 4500, 4668, 4695,
4875, 6111, 6261, 7041, 7320, 7470, 8466, 8652, 8745, 9096, 9345, 9891, 9990, 10050, 10305,
10431, 11196, 13632, 13671, 14766, 15351, 16191, 16341, 16353, 16695, 18480, 18621, 19026,
19566, 22200, 22695, 22956, 27951, 35805, 43560, 44331, 47295, 60030, 90441.

\medskip

 Conjecture 1.1 in the case $a=1$ and $b\in\{0,2\}$ can be refined as follows.

\proclaim{Conjecture 3.2} {\rm (i)} Each natural
number $n>43473$ can be written in the form $2p+T_x$ with $x\in\Z$, where $p$
is zero or a prime.

{\rm (ii)} Any integer $n>636471$ can be written in the form
$2p+T_x$ with $x\in\Z$, where $p$ is zero or a prime congruent to
$1$ modulo $4$. Also, any integer $n>719001$ can be written in the
form $2p+T_x$ with $x\in\Z$, where $p$ is zero or a prime congruent
to $3$ modulo $4$.
\endproclaim

\Remark\ 3.2. We have verified the conjecture for $n\ls 10,000,000$.
As any natural number $n\ls  636,471$ not in the exceptional set $E$
given in Conjecture 1.2 is either a triangular number or a sum of
two odd squares and a triangular number, Conjecture 3.2(ii) implies
the second part of Conjecture 1.2 since any prime $p\eq1\ (\mo\ 4)$
can be written in the form $x^2+y^2$ with $x$ even and $y$ odd.

\medskip

Below is the full list of natural numbers not exceeding 43,473 that cannot be written in the form
$2p+T_x$, where $p$ is 0 or a prime, and $x$ is an integer.
\medskip

2, 8, 18, 30, 33, 57, 60, 99, 108, 138, 180, 183, 192, 240, 243, 318, 321, 360,
366, 402, 421, 429, 495, 525, 546, 585, 591, 606, 693, 696, 738, 831, 840, 850,
855, 900, 912, 945, 963, 1044, 1086, 1113, 1425, 1806, 1968, 2001, 2115, 2190,
2550, 2601, 2910, 3210, 4746, 5013, 5310, 5316, 5475, 5853, 6576, 8580, 9201,
12360, 13335, 16086, 20415, 22785, 43473.

\medskip

For the case $a=2$ and $b=0,2$ of Conjecture 1.1, we have the
following concrete conjecture.

 \proclaim{Conjecture 3.3} {\rm (i)} Any integer $n>849,591$ can be written in the form $4p+T_x$ with $x\in\Z$,
 where $p$ is zero or a prime.

 {\rm (ii)}  Each integer $n>7,718,511$ can be written in the form $4p+T_x$ with $x\in\Z$,
 where $p$ is either zero or a prime congruent to $1$ mod $4$. And each integer $n>6,276,705$
 can be written in the form $4p+T_x$ with $x\in\Z$,
 where $p$ is either zero or a prime congruent to $3$ mod $4$.
 \endproclaim

 \Remark\ 3.3. We have verified Conjecture 3.3 for  $n\ls 30,000,000$.

\medskip

For $a\in\N$  we define  $f(a)$ to be the largest integer not in the form $2^ap+T_x$,
where $p$ is zero or a prime, and
$x$ is an integer. Our conjectures 1.1 and 3.2-3.3, and related computations suggest that
$$f(0)=216,\ f(1)=43473,\ f(2)=849591.$$

\medskip

 Concerning Conjecture 1.4 in the cases $b=2,3$ we have the following concrete conjecture.

\proclaim{Conjecture 3.4} {\rm (i)} Let $n>1$ be an odd integer. Then $n$
can be written in the form $p+x(x+1)$ with $p$ a prime congruent to $1$ mod $4$
and $x$ an integer, if and only if $n$ is not among the following $30$ multiples of three:
$$\align&3,\, 9,\, 21,\, 27,\, 45,\, 51,\, 87,\, 105,\, 135,\, 141,
\\&189,\, 225,\, 273,\, 321,\, 327,\,471,\, 525,\,
627,\,741,\, 861,
\\&975,\,
1197,\, 1461,\, 1557,\, 1785,\, 2151,\, 12285,\, 13575,\, 20997,\, 49755.
\endalign$$
Also, $n$
can be written in the form $p+x(x+1)$ with $p$ a prime congruent to $3$ mod $4$
and $x$ an integer, if and only if $n$ is not among the following $15$ multiples of three:
$$\align&57,\, 111,\, 297,\, 357,\, 429,\, 615,\, 723,\, 765,
\\&1185,\, 1407,\, 2925,\, 3597,\, 4857,\, 5385,\, 5397.
\endalign$$

{\rm (ii) For each $r\in\{1,3,5,7\}$, any odd integer $n>M_r$ can be written
in the form $p+x(x+1)$ with $p$ a prime congruent to $r$ mod $8$ and $x$ an integer,
where $$M_1=358245,\ M_3=172995,\ M_5=359907,\ M_7=444045.$$
\endproclaim

\Remark\ 3.4. We have verified Conjecture 3.4 for odd integers below
$5\times10^6$. It is curious that all the exceptional numbers in the
first part of Conjecture 3.4 are multiples of three.

\medskip

\heading{4. Additional remarks on Conjectures 1.2 and 1.3}\endheading

 As usual, we set
 $$\varphi(q)=\sum_{n=-\infty}^{\infty}q^{n^2}\ \ \t{and}\ \ \psi(q)=\sum_{n=0}^{\infty}q^{T_n}\ \ \ (|q|<1).$$
 There are many known relations between these two theta functions (cf. Berndt [B, pp.\,71-72]). For a $q$-series $F(q)$ we use
 $[q^n]F(q)$ to denote the coefficient of $q^n$ in $F(q)$. By the generating function method,
 Conjecture 1.2 tells that
 $$\gather [q^n]\varphi^2(q^4)\psi(q)>0\ \ \ \t{for any}\ n>864,
 \\ [q^n]\varphi(q^4)\varphi(q^{16})\psi(q)>0\ \ \ \t{for any}\ n>2577,
 \\[q^n](1+q^2\psi^2(q^8))\psi(q)>0\ \ \ \t{for any}\ n>1029.
 \endgather$$

  Here are some of our observations concerning Conjecture 1.3
arising from numerical computations up to $10^6$:

$$S_3=\{4,\, 2578\},\ S_4=\{39\},\ S_{10}=\{87,\, 219,\, 423\},\ S_{60}=\{649,\,1159\};$$
$$S_{15}=\{16,\,49,\,77,\,91,\,136,\,752,\,808,\,931\},\ S_{18}=\{803\};$$
$$S_{24}=\{25,\,49,\,289,\,889,\,1585\},\ S_{36}=\{85,\,91,\,361,\,451,\,1501\};$$
$$\align S_{48}=\{&49,\,125,\,133,\,143,\,169,\,209,\,235,\,265,\,403,\,473,\,
\\&815,\,841,\,1561,\,1679,\,4325,\,8075,\,14953\}.\endalign$$
$$\align S_3^{(1)}=\{&5,\,8,\,11,\,16,\,20,\,50,\,53,\,70,\,113,\,128,\,133,\,200,\,233,
\\&\,245,\,275,\,350,\,515,\,745,\,920,\,1543,\,1865,\,2158,\,3020\}.
\endalign$$
$$S_8(1)=\{1,\,4,\,7,\,16,\,28,\,46,\,88,\,91,\,238,\,373,\,1204\},\ S_8(5)=\{26,\,65,\,176\};$$
$$S_9(1)=\{1,\,6,\,16,\,141\},\ S_9(4)=\{5,\, 19,\,50,\,75\}, \ S_9(7)=\{2,\,73,\,98,\,232,\,448\}.$$

\bigskip

\noindent{\bf Added in proof}. The author has set up a webpage
devoted to mixed sums of primes and other terms with the website

\qquad\qquad {\tt http://math.nju.edu.cn/$\sim$zwsun/MSPT.htm}.

\enddocument